\documentclass[onecolumn,a4paper]{revtex4}
\usepackage[utf8]{inputenc}
\usepackage[english]{babel}

\usepackage{amsmath}    % need for subequations
\usepackage{verbatim}   % useful for program listings
\usepackage{hyperref}   % use for hypertext links, including those to external documents and URLs
\usepackage{bbold}
\usepackage{mathrsfs}
\usepackage{textcomp}
\usepackage{amsthm}
\usepackage{amssymb}
\usepackage[capitalize]{cleveref}
\allowdisplaybreaks

\usepackage{float}

\usepackage{tikz}
\usetikzlibrary{calc,decorations.markings}

\usepackage{scrextend}
\changefontsizes[15pt]{11pt}

\begin{document}
\title{Integrals containing the infinite product $\prod_{n=0}^\infty\left[1+\left(\frac{x}{b+n}\right)^3\right]$} 
\author{Martin Nicholson}
\begin{abstract}  
      We study several integrals that contain the infinite product ${\displaystyle\prod_{n=0}^\infty}\left[1+\left(\frac{x}{b+n}\right)^3\right]$ in the denominator of their integrand. These considerations lead to closed form evaluation $\displaystyle\int_{-\infty}^\infty\frac{dx}{\left(e^x+e^{-x}+e^{ix\sqrt{3}}\right)^2}=\frac{1}{3}$ and to some other formulas.
\end{abstract}
\maketitle
 
$\mathbf{1}.$ The infinite product
$$
\prod_{n=0}^\infty\left[1+\left(\frac{\alpha+\beta}{n+\alpha}\right)^3\right]
$$
and more general products have been studied in the literature (see [\onlinecite{berndt}], ch. 16). In this paper we consider integrals of the form
\begin{align}\label{int}
\int_{0}^\infty P_b(x)f(x)dx,
\end{align}
where 
\begin{equation}\label{product}
P_b(x)=\frac{1}{\prod\limits_{k=0}^\infty\left(1+\frac{x^3}{(k+b)^3}\right)}.
\end{equation}
Several such integrals will be evaluated in closed form. However while others do not have a closed form will allow us to evaluate some integrals of elementary functions.

Note that the infinite product in \eqref{product} can be written in terms of Gamma functions [\onlinecite{ww}]
$$
P_b(x)=\frac{\Gamma(b+x)\Gamma(b+\omega x)\Gamma(b+x/\omega)}{\Gamma^3(b)},\quad \omega=e^{\frac{2\pi i}{3}}.
$$
The notation $\omega=e^{\frac{2\pi i}{3}}$ for third root of unity will be used throughout the paper.

$\mathbf{2}.$ Consider the contour integral
\begin{align}
\int_C P_b(z)\frac{dz}{z}.
\end{align}
along the contour depicted in Fig.1. We assume that $b>0$. The most interesting cases considered in this paper correspond to $b=1$ and $b=1/2$.
\begin{figure}
\centering
    \begin{tikzpicture}
  [
    decoration={%
      markings,
      mark=at position 1.6cm with {\arrow[line width=1pt]{>}},
      mark=at position 0.45 with {\arrow[line width=1pt]{>}},
      mark=at position 0.83 with {\arrow[line width=1pt]{>}},
      mark=at position -5mm with {\arrow[line width=1pt]{>}},
    }
  ]
  \draw [help lines,->] (-2.2,0) -- (4,0) coordinate (xaxis);
  \draw [help lines,->] (0,-0.5) -- (0,4) coordinate (yaxis);
  \path [draw, line width=0.8pt, postaction=decorate] (0.5,0) node [below] {$\varepsilon$} -- (3.5,0) node [below] {$R$} arc (0:120:3.5) node [left] {$R\omega$} -- (-.25,0.43) node [left] {$\varepsilon \omega$} arc (120:0:.5);
  \node [below left] {$O$};
  \node at (0.45,.7) {$C_{\varepsilon}$};
  \node at (2,3.3) {$\Gamma_{R}$};
  \node at (0,-0.8) {$\text{Fig.1}$};
\end{tikzpicture}
\end{figure}

Inside the contour of integration, the integrand $h(z)=P_b(z)/z$ has simple poles at $z=-(k+b-1)/\omega,~k\in\mathbb{N}$,  with residues
$$
\frac{(-1)^{k}}{(k-1)!}\frac{|\Gamma(b-\omega(k+b-1))|^2}{(k+b-1)\Gamma^3(b)},
$$
and no poles on the contour of integration if we choose $R=N+b-1/2$ for some large natural number $N$. Also $h(z)dz$ is symmetric under the change $z\to\omega z$, and as a consequence the integrals along straight lines cancel each other out. Let's denote the integrals along $\Gamma_R$ and $C_\varepsilon$ as $I_R$ and $I_\varepsilon$ respectively.
Then
$$
\lim_{\varepsilon\to 0}I_\varepsilon=-\frac{2\pi i}{3},
$$
and (Appendix \ref{a})
$$
\lim_{R\to +\infty}I_R=0.
$$
Using residue theorem we get
\begin{equation}\label{sum}
    \sum_{n=0}^\infty \frac{(-1)^n}{n!}\frac{\left|\Gamma(b-\omega(n+b))\right|^2}{n+b}=\frac13\Gamma^3(b).
\end{equation}
The integral  3.985.1 from [\onlinecite{GR}]
\begin{equation}\label{gr}
    \int_{-\infty}^\infty \frac{e^{iax} \space dx}{\cosh^\nu \beta x} = \frac{2^{\nu- 1}}{\beta\Gamma(\nu)} \Gamma\left( \frac{\nu}{2} + \frac{ai}{2\beta} \right) \Gamma\left( \frac{\nu}{2} - \frac{ai}{2\beta} \right)
\end{equation}
allows to write \eqref{sum} as an integral of a hypergeometric function
\begin{equation}\label{hypergeometric}
    \int_{-\infty}^\infty\frac{e^{ib\sqrt{3}x}}{\cosh^{3b} x}{}_2F_1\left({b,3b \atop b+1}\Biggr|-\frac{e^{i\sqrt{3}x}}{2\cosh x}\right)dx=2^{3b-1}\frac{b}{3}\frac{\Gamma^3(b)}{\Gamma(3b)}.
\end{equation}
$\mathbf{3.}$ Here we specialize $b$ in \eqref{hypergeometric} so that the hypergeometirc function can be written in terms of elementary functions. This happens when $b=1+3n$ or $b=1/2+3n$, where $n$ is a non-negative integer. Only the two cases with $n=0$ are considered below:

Let $b=1/2$, then the hypergeometric function becomes $\sqrt{\frac{{2\cosh x}}{{{2\cosh x}+{e^{i\sqrt{3}x}}}}}$ and we get
\begin{equation}\label{sqrt}
    \int_{-\infty}^{\infty}\frac{\text{sech}x~e^{\frac{i\sqrt{3}}{2}x}dx}{\sqrt{e^x+e^{-x}+e^{i\sqrt{3}x}}}=\frac{\pi}{3}.
\end{equation}

If $b=1$ then the hypergeometric function becomes $\frac{2\cosh x\left(4\cosh x+e^{i\sqrt{3}x}\right)}{\left(2\cosh x+e^{i\sqrt{3}x}\right)^2}$ and we get
$$
    \int_{-\infty}^{\infty}\frac{e^{i\sqrt{3}x}}{\cosh^{2} x}\frac{4\cosh x+e^{i\sqrt{3}x}}{\left(2\cosh x+e^{i\sqrt{3}x}\right)^2}dx=\frac{2}{3},
$$
which due to $\frac{4\cosh x+e^{i\sqrt{3}x}}{\cosh^2x\left(2\cosh x+e^{i\sqrt{3}x}\right)^2}=-\frac{4}{\left(2\cosh x+e^{i\sqrt{3}x}\right)^2}+\frac{1}{\cosh^2 x}$ can be simplified further as
\begin{equation}\label{int1}
    \int_{-\infty}^\infty\frac{dx}{\left(e^x+e^{-x}+e^{ix\sqrt{3}}\right)^2}=\frac{1}{3}.
\end{equation}
It is interesting to note that there is another way to write the sum \eqref{sum} with $b=1$ as an integral
\begin{equation}\label{int2}
    \int_{-\infty}^{\infty}\frac{e^{i \sqrt{3} x} \cosh x}{\left(e^x+e^{-x}+e^{i \sqrt{3} x}\right)^2}dx=\frac{1}{12}.
\end{equation}
One might observe how the $2\pi /3$ rotation symmetry of the product $\prod\limits_{k=1}^\infty\left(1+\frac{x^3}{k^3}\right)$ manifests itself in \eqref{int1} and \eqref{int2}: The set of roots of the equation $e^x+e^{-x}+e^{i \sqrt{3} x}=0$ has the same $2\pi /3$ rotation symmetry (see Appendix \ref{b}).

$\mathbf{4}.$ The last integral in section $3$ gives
\begin{align}
\nonumber \int_{-\infty}^\infty\frac{e^{i\sqrt{3}x}\sinh x \,dx}{\left(2\cosh x+e^{i\sqrt{3}x}\right)^2}&=\frac12\int_{-\infty}^\infty\frac{e^{i\sqrt{3}x}\,d(2\cosh x+e^{i\sqrt{3}x})}{\left(2\cosh x+e^{i\sqrt{3}x}\right)^2}-\frac{i\sqrt{3}}{2}\int_{-\infty}^\infty\frac{e^{2i\sqrt{3}x}dx}{\left(2\cosh x+e^{i\sqrt{3}x}\right)^2}\\
\nonumber &=\frac{i\sqrt{3}}{2}\left(\int_{-\infty}^\infty\frac{e^{i\sqrt{3}x}dx}{2\cosh x+e^{i\sqrt{3}x}}-\int_{-\infty}^\infty\frac{e^{2i\sqrt{3}x}dx}{\left(2\cosh x+e^{i\sqrt{3}x}\right)^2}\right)\\
\nonumber &={i\sqrt{3}}\int_{-\infty}^{\infty}\frac{e^{i \sqrt{3} x} \cosh x}{\left(e^x+e^{-x}+e^{i \sqrt{3} x}\right)^2}dx=\frac{i\sqrt{3}}{12}.
\end{align}
Thus we have
\begin{equation}\label{sinh}
   \int_{-\infty}^\infty\frac{e^{i\sqrt{3}x}\sinh x \,dx}{\left(2\cosh x+e^{i\sqrt{3}x}\right)^2} =\frac{i\sqrt{3}}{12}.
\end{equation}

$\mathbf{5}.$ After the substitution $t=e^{2x}$ equation \eqref{int1} becomes
\begin{equation}
    \int_0^\infty\frac{dt}{(1+t+t^{\,\alpha})^2}=\frac23, \quad \alpha=\frac{1+i\sqrt3}{2},
\end{equation}
while combining \eqref{int2} and \eqref{sinh} we find two analogous representations
\begin{equation}\label{int3}
    \int_0^\infty\frac{t^{\alpha}dt}{(1+t+t^{\,\alpha})^2}=\frac{\alpha}{3},
\end{equation}
\begin{equation}\label{int4}
    \int_0^\infty\frac{t^{\alpha-1}dt}{(1+t+t^{\,\alpha})^2}=\frac{1}{3\alpha}.
\end{equation}
\eqref{int4} is related to \eqref{int3} by complex conjugation and change of variable $t\to 1/t$.

$\mathbf{6}.$ If we apply the approach of section 2 to
\begin{align*}
\int_C P_b(z)dz,
\end{align*}
instead, then the integrals over straight lines no longer cancel out. However, there is nevertheless a simplification: in this case the sum analogous to \eqref{sum} reduces to an integral of elementary function for all $b$ so that in this case we find the transformation
\begin{align}\label{transformation}
\int_{0}^\infty\frac{dx}{\prod\limits_{k=0}^\infty\left(1+\frac{x^3}{(k+b)^3}\right)}=\frac{4\pi\Gamma(3b)}{\Gamma^3(b)\sqrt{3}}\int_{-\infty}^\infty\frac{e^{ixb\sqrt{3}}\ dx}{\left(e^x+e^{-x}+e^{ix\sqrt{3}}\right)^{3b}}.
\end{align}

$\mathbf{7}.$ In principle integrals \eqref{sqrt} and \eqref{int1} can be written in terms of real-valued functions by calculating the real part of the integrand. The resulting formulas are cumbersome and therefore omitted. However there is another way to get a compact integral of a real valued function, at least for \eqref{int1}. First of all, all the roots of the function $~e^z+e^{-z}+e^{iz\sqrt{3}}~$  lie on the three rays $z=ir\omega^k,~r>0$, $(k=0,1,2)$ (see Appendix \ref{a}). If one bends the contour of integration so that it never crosses these zeroes then the integral \eqref{int1} will not change. Since the integrand decreases exponentially when $z\to \infty$, $0<\arg z<\pi/6$ or $5\pi/6<\arg z<\pi$ we have
$$
\frac{1}{\beta}\int_{0}^\infty\frac{dx}{\left(e^{-x/\beta}+e^{x/\beta}+e^{-ix\sqrt{3}/\beta}\right)^2}+\beta\int_{0}^\infty\frac{dx}{\left(e^{\beta x}+e^{-\beta x}+e^{i\beta x\sqrt{3}}\right)^2}=\frac{1}{3},\quad \beta=e^{\pi i/6},
$$
and after elementary simplifications
\begin{equation}
    \int_0^\infty\frac{e^{x\sqrt{3}}\cos\left(\frac\pi6-x\right)}{\left(2\cos x+e^{x\sqrt{3}}\right)^2}dx=\frac16.
\end{equation}
Similarly, for the case $b=1$ of \eqref{transformation}
\begin{equation}
     \int_{0}^\infty\frac{dx}{\left(1+\frac{x^3}{1^3}\right)\left(1+\frac{x^3}{2^3}\right)\left(1+\frac{x^3}{3^3}\right)\ldots}=8\pi\int_0^\infty\frac{e^{x\sqrt{3}}~dx}{\left(2\cos x+e^{x\sqrt{3}}\right)^3}.
 \end{equation}
 
$\mathbf{8}.$ It turns out that \eqref{int1} has a parametric extention. Consider the contour integral
\begin{align}\label{ci}
\int_{C'}P_1(z)\frac{e^{a z}dz}{z},
\end{align}
where the contour $C'$ is a circle of radius $R=N+1/2$ for large natural $N$. Since $|P_1(z)|$ decreases exponentially with $N$ on the circle $C'$ (Appendix \ref{a}), \eqref{ci} will be zero in the limit $N\to\infty$ for sufficiently small $|a|$. Therefore the sum of residues of the integrand over three sets of simple poles $z=-ke^{2\pi ij/3}$,$~k\in\mathbb{N}$, $(j=0,1,2)$ plus a simple pole at the origin, will be $0$ according to residue theorem. As a result one will obtain three sums similar to \eqref{sum} and then convert them to integrals of the type \eqref{int1}. However there is a trick that allows to avoid these calculations. Note that the factor $e^{az}$ in \eqref{ci} will introduce additional factors $\exp{\left(-ake^{\frac{2\pi i j}{3}}\right)}$ in the sum over residues. When converted to an integral via \eqref{gr} these factors have the effect of multiplying $e^{ix\sqrt{3}}$ by $\exp{\left(-ae^{\frac{2\pi i j}{3}}\right)}$:
\begin{align*}
 \sum_{j=1}^3\int_{-\infty}^\infty\frac{dx}{\left(e^x+e^{-x}+\exp{\left(-ae^{\frac{2\pi i j}{3}}\right)}e^{ix\sqrt{3}}\right)^2}=1,
\end{align*}
or equivalently
\begin{align}
\int\limits_{-\infty}^\infty\frac{dx}{\left(e^x+e^{-x}+e^{a+ix\sqrt{3}}\right)^2}+\int\limits_{-\infty}^\infty\frac{e^adx}{\left(e^{a+x}+e^{-x}+e^{ix\sqrt{3}}\right)^2}+\int\limits_{-\infty}^\infty\frac{e^adx}{\left(e^{a+x}+e^{-x}+e^{-ix\sqrt{3}}\right)^2}=1,
\end{align}
where $|a|$ is sufficiently small.

$\mathbf{9}.$ There is a similarity between \eqref{sum} with $b=1$ and the identity due to Ramanujan ([\onlinecite{berndt1}], p. 309)
$$
e^{ay}=\frac{-a}{2ci}\sum_{k=0}^\infty\frac{\Gamma\left(\frac{-a+k(ci-b)}{2ci}\right)(-2ie^{-by}\sin cy)^k}{\Gamma\left(\frac{-a-k(ci-b)}{2ci}+1\right)k!}.
$$
After equating the coefficents of $a^1$ in Taylor series expansion of both sides into powers of $a$ and transforming the Gamma function in the denominator via Euler's reflection formula
\begin{equation}\label{r}
    y=\frac{1}{2\pi ic}\sum_{k=1}^\infty\frac{(-1)^k}{k!}\Gamma\left(\frac{k}{2}-\frac{ikb}{2c}\right)\Gamma\left(\frac{k}{2}+\frac{ikb}{2c}\right)\sin\pi\left(\frac{k}{2}+\frac{ikb}{2c}\right)(-2ie^{-by}\sin cy)^k.
\end{equation}
To make this similarity more exact we differentiate \eqref{r} with respect to $y$, divide by $(c\cot cy-b)$, repeat this procedure one more time and then set $c=1$, $b=\sqrt{3}$
\begin{equation}\label{r2}
    \frac{2\pi i\sin y}{(\cos y-\sqrt{3}\sin y)^3}=\sum_{k=1}^\infty\frac{(-1)^k}{k!}{|{\Gamma\left(1-\omega k\right)}|^2}\cdot\sin\left(\pi ke^{\pi i/3}\right)(-2ie^{-\sqrt{3}y}\sin y)^k.
\end{equation}
This series converges when $|2e^{-\sqrt{3}y}\sin y|<e^{-\frac{\pi}{2\sqrt{3}}}$ [\onlinecite{berndt1}]. It will be convenient to use another variable $\alpha$ related to $y$ by
\begin{equation*}
e^{-\alpha}=2e^{-\sqrt{3}y}\sin y.
\end{equation*}
The condition that the series \eqref{r2} converges now takes a very simple form $~\text{Re}\phantom{.}\alpha>\frac{\pi}{2\sqrt{3}}$. In the following it will be assumed for simplicity that $\alpha>\frac{\pi}{2\sqrt{3}}$.

Is it possible that \eqref{r2} leads to evaluation of integrals with infinite product $\prod\limits_{k=1}^\infty\left(1+\frac{z^3}{k^3}\right)$? Consider the contour integral
\begin{align}
\int\limits_C\frac{\left(-ie^{-\alpha}\right)^{-\omega z}\sin\pi z }{z\prod\limits_{k=1}^\infty\left(1+\frac{z^3}{k^3}\right)}~dz,
\end{align}
where $C$ is the contour in Fig.2. Due to the asymptotics 
$$
|\left(-ie^{-\alpha}\right)^{-\omega z}\sin\pi z| \sim \frac12 \exp\left[-\left(\frac{\alpha}{2}+\frac{\pi\sqrt{3}}{4}\right)x+\left(\frac{5\pi}{4}-\frac{\alpha\sqrt{3}}{2}\right)y\right],\quad 0<\arg z<\frac{2\pi}{3},
$$
and the result of Appendix \ref{a}, the integral over circular arc $\Gamma_R$ vanishes in the limit $R\to\infty$ if
$$
\begin{cases}
-\left(\frac{\alpha}{2}+\frac{11\pi}{4\sqrt{3}}\right)x+\left(\frac{5\pi}{4}-\frac{\alpha\sqrt{3}}{2}\right)y<0 ,\quad 0<\arg z<\frac{\pi}{3}, \\
\left(\frac{\pi}{2\sqrt{3}}-\alpha\right)(x+y\sqrt3)<0,\phantom{..}~~\qquad \qquad \frac{\pi}{3}<\arg z<\frac{2\pi}{3}.
\end{cases}
$$
When $~\alpha>\frac{\pi}{2\sqrt{3}}~$ these conditions are satisfied automatically. 

The same approach as in section $2$ yields
\begin{align*}
    &\int_0^\infty\frac{\left(-i e^{-\alpha}\right)^{-\omega x}\sin\pi x -\left(-ie^{-\alpha}\right)^{- x/\omega}\sin\pi \omega x}{x\prod\limits_{k=1}^\infty\left(1+\frac{x^3}{k^3}\right)}~dx\\
    &=2\pi i\sum_{k=1}^\infty\frac{(-1)^{k-1}}{k!}{|{\Gamma\left(1-\omega k\right)}|^2}\cdot\sin\left(\pi ke^{\pi i/3}\right)(-ie^{-\alpha})^k\\
    &=\frac{4\pi^2\sin y}{(\cos y-\sqrt{3}\sin y)^3}.
\end{align*}
For real $\alpha$ one can decompose the function in the numerator of this integral into real and imaginary parts
$$
\left(-i e^{-\alpha}\right)^{-\omega x}\sin\pi x -\left(-ie^{-\alpha}\right)^{- x/\omega}\sin\pi \omega x=f(x,\alpha)+ig(x,\alpha)
$$
where
\begin{equation}
    \nonumber f(x,\alpha)=\frac12 e^{-{ \sqrt{3} \pi  x}/{4}-\alpha x/2}\left(e^{\sqrt{3} \pi  x} \sin\frac{\left(\pi -2 \sqrt{3} \alpha\right)x}{4} +2 \sin\frac{\left(2 \sqrt{3} \alpha+3 \pi \right) x}{4}-\sin\frac{\left(2 \sqrt{3} \alpha-5 \pi \right) x}{4}\right),
\end{equation}
\begin{equation}
    \nonumber g(x,\alpha)=\frac12 e^{-{ \sqrt{3} \pi  x}/{4}-\alpha x/2}\left(\cos\frac{ \left(2 \sqrt{3} \alpha-5 \pi \right) x}{4}-e^{\sqrt{3} \pi  x} \cos\frac{\left(\pi -2 \sqrt{3} \alpha\right) x}{4}\right).
\end{equation}
As a result
\begin{equation}
    \int_{0}^\infty\frac{g(x,\alpha)\ dx}{x\prod\limits_{k=1}^\infty\left(1+\frac{z^3}{k^3}\right)}=0,
\end{equation}
\begin{equation}
    \int_{0}^\infty\frac{f(x,\alpha)\ dx}{x\prod\limits_{k=1}^\infty\left(1+\frac{z^3}{k^3}\right)}=\frac{8\pi^2\sin y}{\left(\sqrt{3}\sin y-\cos y\right)^3},
\end{equation}
 where $y$ is the root of the equation $2e^{-y\sqrt{3}}\sin y=e^{-\alpha}$ near $y=0$.
 
These formulas simplify when $\alpha =\frac{5\pi}{2\sqrt{3}}$
\begin{equation}
    \int\limits_{0}^\infty\frac{\left(1-e^{\pi\sqrt{3}x}\cos\pi x\right)e^{-\frac{2\pi}{\sqrt{3}}x}\ dx}{x\left(1+\frac{x^3}{1^3}\right)\left(1+\frac{x^3}{2^3}\right)\left(1+\frac{x^3}{3^3}\right)\ldots}=0,
\end{equation}
\begin{equation}
    \int\limits_{0}^\infty\frac{\sin\pi x\left(4\cos\pi x-e^{\pi\sqrt{3}x}\right)e^{-\frac{2\pi}{\sqrt{3}}x}dx}{x\left(1+\frac{x^3}{1^3}\right)\left(1+\frac{x^3}{2^3}\right)\left(1+\frac{x^3}{3^3}\right)\ldots}=\frac{8\pi^2\sin y}{\left(\sqrt{3}\sin y-\cos y\right)^3},
\end{equation}
 where $y=0.0054167536\ldots$ is the root of the equation $2e^{-y\sqrt{3}}\sin y=e^{-\frac{5\pi}{2\sqrt{3}}}$.
 
 $\mathbf{10}.$ The hyperbolic log-trigonometric integral
 \begin{equation}\label{logtrig}
     \text{Im}\int_0^\infty \frac{dt}{\left(i t\sqrt{3}+\ln(2\sinh t)\right)^2}=0,
 \end{equation}
 or in terms of real valued functions
 \begin{equation}\label{real}
     \int_{0}^{\infty}\frac{t\ln\left(2\sinh t\right)}{\left[3t^{2} + \ln^{2}\left(2\sinh t\right)\right]^{2}}\,{d}t = 0,
 \end{equation}
is also related to the infinite product in the title. Indeed
\begin{align}\label{omega}
\nonumber\int_{0}^\infty\frac{\sin\pi x~dx}{x\prod\limits_{k=1}^\infty\left(1-\frac{x^3}{k^3}\right)}
\nonumber &=\int_{0}^\infty \frac{\sin\pi x}{x}\Gamma(1-x)\Gamma(1-\omega x)\Gamma(1-x/\omega)dx\\
\nonumber &=\pi\int_{0}^\infty \frac{ \Gamma(1-\omega x)\Gamma(1-x/\omega)}{\Gamma(1+x)}dx\\
\nonumber &=\pi\int_{0}^\infty B(1-\omega x,1-x/\omega)dx\\
&=\pi\int_{0}^\infty dx\int_0^\infty \frac{t^{-\omega x}}{(1+t)^{x+1}}dt
\end{align}
Changing the order of integration and calculating the integral over $x$ we get
\begin{equation}\label{int5}
\int_{0}^\infty\frac{\sin\pi x~dx}{x\prod\limits_{k=1}^\infty\left(1-\frac{x^3}{k^3}\right)}
=-2\pi \int_0^\infty \frac{dt}{\left(i t\sqrt{3}+\ln(2\sinh t)\right)^2}.
\end{equation}
\eqref{logtrig} is the statement of the fact that the integral on the RHS of \eqref{int5} is real. Of course by replacing in \eqref{omega} $\omega$ with any complex number of unit argument one gets other integrals like \eqref{logtrig}.

It is known that Laplace transform of the digamma function leads to some log-triginometric integrals [\onlinecite{glasser,moll,dixit}] that contain the expression
$x^2+\ln^2(2e^{-a}\cos x)$ in the denominator. This expression should be compared to the expression $3t^{2} + \ln^{2}\left(2\sinh t\right)$ in the denominator of \eqref{real}.

\appendix\section{Asymptotics of the product of gamma functions}\label{a}
Due to the asymptotic relation
$$
\ln\Gamma(z)=\left(z-\frac12\right)\ln z-z+O(1),\quad |\arg z|<\pi,
$$
one has
$$
\ln\left\{\Gamma(b+z)\Gamma(b+\omega z)\Gamma(b+z/\omega)\right\}=3\left(b-\frac12\right)\ln z-\frac{2\pi}{\sqrt{3}}z+O(1),\quad |\arg z|<\frac{\pi}{3}.
$$
From this it follows that
$$
|P_b(z)|=C|z|^{3b-3/2}\cdot\begin{cases}
e^{-\frac{2\pi}{\sqrt{3}}x},~0<\arg z<\frac{\pi}{3},\\
e^{\tfrac{\pi}{\sqrt{3}}x-\pi y},~\frac{\pi}{3}<\arg z<\frac{2\pi}{3},
\end{cases}
$$
where $z=x+iy$.

\section{Roots of the equation $e^{i\sqrt3 z}+2\cosh z=0$}\label{b}

The fact that the roots of the equation $e^{i\sqrt3 z}+2\cosh z=0$ are symmetric under $z\to \omega z$ is easy to check directly.

Since $\frac12 e^{-\pi\sqrt3 /2}=0.0329...$ is quite small the equation $e^{i\sqrt3 z}+2\cosh z=0$ will have roots close to $\pi i \left(n+\frac12\right)$, where $n$ is a non-negative integer. Below it is shown that these are the only roots in the upper half plane.

Let $f(z)=e^{i\sqrt3 z}$,$~g(z)=2\cosh z$. Obviously, on the real axis $|f(z)|<|g(z)|$. Now consider $f(z)$ and $g(z)$ on the closed contour $C$ depicted in Fig.\hyperref[figure2]{2}
\begin{figure}[H]
  \centering
  \begin{tikzpicture}
  \draw [help lines,->] (-4,0) -- (4,0) coordinate (xaxis);
  \draw [help lines,->] (0,-0.5) -- (0,4) coordinate (yaxis);
  \path [draw, line width=0.8pt, postaction=decorate] (0,0) -- (3.5,0) node [below] {$R$} arc (0:180:3.5) node [below] {$-R$} -- (0,0) ;
  \node [below left] {$O$};
  \node at (2,3.3) {$\Gamma_{R}$};
  \node at (0,-0.8) {$\text{Fig.2}$};
\end{tikzpicture}
\end{figure}
Here $\Gamma_R$ is a semicircle of radius $R=\pi N$ for some large natural number $N$.
We have for $z=x+iy\in \Gamma_R$
$$
|f(z)|=e^{-\sqrt3 y}\le 1,
$$
$$
|g(z)|=2\sqrt{\sinh^2x+\cos^2\sqrt{\pi^2N^2-x^2}}\ge 2.
$$
Thus $|f(z)|<|g(z)|$ on the contour $C$. According to Rouche's theorem this means that the function $f(z)+g(z)$ has the same number of roots inside the contour $C$ as the function $g(z)$, as required.

This analysis shows that the roots of $e^{i\sqrt3 z}+2\cosh z=0$ are located on the three rays $z=ir\omega^k$, $~k=0,1,2$.

\end{document}